\documentclass[10pt]{amsart}
\usepackage{amsmath,amsthm}
\usepackage[colorlinks=true,linkcolor=blue,urlcolor=blue,citecolor=blue]{hyperref}
\usepackage[all]{xy}

\begin{document}

\newtheorem{theorem}{Theorem}[section]
\newtheorem{lemma}[theorem]{Lemma}
\newtheorem{proposition}[theorem]{Proposition}
\newtheorem{corollary}[theorem]{Corollary}

\theoremstyle{definition}
\newtheorem{definition}[theorem]{Definition}
\newtheorem{remark}[theorem]{Remark}
\newtheorem{example}[theorem]{Example}
\newtheorem{problem}[theorem]{Problem}
\newtheorem*{acknowledgements}{Acknowledgements}

\numberwithin{equation}{section}

\renewcommand{\subjclassname}{%
\textup{2020} Mathematics Subject Classification}

\numberwithin{equation}{section}  

\newtheorem{conjecture}[theorem]{Conjecture}


\title{Arithmetic Properties For $(r,s)$-Regular Partition Functions With Distinct Parts}

\author[R. Drema and N. Saikia]{Rinchin Drema and Nipen Saikia$^\ast$}
\address{Department of Mathematics, Rajiv Gandhi University, Rono Hills, Doimukh, Arunachal Pradesh, India, Pin-79111.}
\email{dremarinchin3@gmail.com; nipennak@yahoo.com}

\keywords{regular partition, distinct parts, $q$-series identities, Ramanujan's theta-functions.\\
	            $\ast$ \textit{Corresponding author}.}

 \maketitle

\begin{abstract}
For any relatively prime integers $r$ and $s$, let $a_{r,s}(n)$ denote the number of $(r,s)$-regular partitions of a positive integer of $n$ into distinct parts. Prasad and Prasad (2018) proved many infinite families of congruences modulo 2 for $a_{3,5}(n)$. In this paper, we establish families of congruences modulo 2 and 4 for $a_{r,s}(n)$ with $(r,s)\in$ \{(2,5), (2,7), (4,5), (4,9)\}.  For example, we show that for all  $\beta \geq 0$ and  $n \geq 0,$ we have
$$a_{2,5}\Big(4\cdot 5^{2\beta+1}n+\dfrac{37\cdot5^{2\beta}-1}{6}\Big)\equiv0 \pmod4. $$ 
\end{abstract}
\section{Introduction} The partition of a positive integer $n$ is a non-increasing sequence of positive integers whose sum is equal to $n$. If  $p(n)$ denote the number of partitions of a positive integer $n$ and we adopt the convention $p(0)=1$, then the generating function for $p(n)$ satisfies the identity
\begin{equation}\label{pngeni}\sum_{n=0}^{\infty}p(n)q^n=\dfrac{1}{(q;q)_{\infty}},\end{equation}where   
\begin{equation}\label{9}
(a;q)_{\infty}=\prod_{n=0}^{\infty}(1-aq^n).
\end{equation} 
Throughout this paper, we write  
$$f_k:=(q^k;q^k)_{\infty}, \quad \text{for any integer}\quad k\ge 1.$$
Ramanujan \cite{rd} established the following beautiful  congruences for all $n\geq0$:
$$p(5n+4) \equiv0\pmod5,\quad p(7n+5)\equiv0 \pmod7\quad  \text{and}\quad 
p(11n+6)\equiv0\pmod {11}.$$
Ramanujan's congruences on $p(n)$ have motivated many mathematicians to seek similar results for restricted partition functions. One example is the $\ell$-regular partition function $b_{\ell}(n)$, which counts the number of partitions of $n$ in which no part is divisible by $\ell$ and whose generating function satisfies the identity
\begin{equation}\label{lregularintro}
\sum_{n=0}^{\infty} b_{\ell}(n)q^n=\frac{f_\ell}{f_1}.\end{equation}
 Many results on the arithmetic of $b_{\ell}(n)$ have been established (see, for example, \cite{Cal33,P34,W44}). 
 
 For relatively prime integers $r$ and $s$, an $(r,s)$-regular partition is one in which  none of the parts is divisible by $r$ or $s$. Denote by $a_{r,s}(n)$, the number of $(r,s)$-regular partitions of $n$ into distinct parts. For example, $a_{2,5}(13)=2$ since the $(2,5)$-regular partitions of 13 into distinct parts are $13$ and $9+3+1$. The generating function for $a_{r,s}(n)$ satisfies the identity  
 \begin{equation}\label{eq1}
 \sum_{n=0}^{\infty}a_{r,s}(n)q^n=\dfrac{(-q;q)_{\infty}(-q^{rs};q^{rs})_{\infty}}{(-q^r;q^r)_{\infty}(-q^s;q^s)_{\infty}}.\end{equation}  
 Prasad and Prasad \cite{MK} proved many infinite families of congruences modulo 2 for $a_{3,5}(n)$.

 In this paper, we establish families of congruences modulo 2 for $a_{2,5}(n)$, $a_{2,7}(n)$, $a_{4,5}(n)$ and $a_{4,9}(n)$. We also prove congruences modulo 4 for $a_{2,5}(n)$. The congruences are listed in the following theorems: 
 
\begin{theorem}\label{thm1}
For every $n\geq0$, we have 
\begin{equation}\label{e1}
a_{2,5}(4n+2)\equiv0 \pmod2
\end{equation} and 
\begin{equation}\label{e2}
a_{2,5}(4n)
\equiv \left\{
        \begin{array}{ll}
         1 \pmod2,   
           & if~ n~ is~ a~ pentagonal~ number \\
           0 \pmod2,
           & otherwise. 
        \end{array}
    \right.
\end{equation}
\end{theorem}

\begin{theorem} \label{thm3}
Let $p>5$ be a prime with $\Big(\dfrac{-10}{p}\Big)=-1$ and $1\leq j\leq p-1$. Then for all $\gamma\ge 0$, we have
\begin{equation}\label{e8}
\sum _{n=0}^{\infty} a_{2,5}\Big(4\cdot p^{2\gamma}n+\dfrac{7\cdot p^{2\gamma}-1}{6}\Big)q^n\equiv f_2f_5 \pmod2,
\end{equation}
\begin{equation}\label{e9}
 a_{2,5}\Big(4\cdot p^{2\gamma+1}(pn+j)+\dfrac{7\cdot p^{2\gamma+2}-1}{6}\Big)\equiv 0 \pmod2,
\end{equation}
\begin{equation}\label{e8a}
\sum _{n=0}^{\infty} a_{2,5}\Big(20\cdot p^{2\gamma}n+\dfrac{55\cdot p^{2\gamma}-1}{6}\Big)q^n\equiv f_1f_{10} \pmod2
\end{equation} and
\begin{equation}\label{e9b}
 a_{2,5}\Big(20\cdot p^{2\gamma+1}(pn+j)+\dfrac{55\cdot p^{2\gamma+2}-1}{6}\Big)\equiv 0 \pmod2.
\end{equation}
\end{theorem}

\begin{theorem}\label{thm4}
If $w_1\in{ \left\lbrace 13, 37\right\rbrace}$~and~$w_3\in{ \left\lbrace 41,89\right\rbrace}$. Then for all $\beta\geq0$, we have
\begin{equation}\label{e10}
\sum _{n=0}^{\infty} a_{2,5}\Big(4\cdot 5^{2\beta}n+\dfrac{13\cdot 5^{2\beta}-1}{6}\Big)q^n\equiv2q^2 f_1f_{20}^3 \pmod4
,\end{equation}
\begin{equation}\label{e11}
\sum _{n=0}^{\infty} a_{2,5}\Big(4\cdot 5^{2\beta+1}n+\dfrac{17\cdot 5^{2\beta+1}-1}{6}\Big)q^n\equiv2 f_4^3 f_5\pmod4,
\end{equation}
\begin{equation}\label{e12}
 a_{2,5}\Big(4\cdot 5^{2\beta+1}n+\dfrac{w_1\cdot 5^{2\beta}-1}{6}\Big)\equiv0 \pmod4
\end{equation} and
\begin{equation}\label{e13}
a_{2,5}\Big(4\cdot 5^{2(\beta+1)}n+\dfrac{w_3\cdot 5^{2\beta+1}-1}{6}\Big)\equiv0 \pmod4.
\end{equation}
\end{theorem}

\begin{theorem}\label{thm5} Let $p>7$ be a prime with $\Big(\dfrac{-14}{p}\Big)=-1$ and $1\leq j\leq p-1$. Then for all $\alpha\ge 0$, we have
\begin{equation}\label{f1}
\sum _{n=0}^{\infty} a_{2,7}\Big(2\cdot p^{2\alpha}n+\dfrac{5\cdot p^{2\alpha}-1}{4}\Big)q^n\equiv f_1f_{14}\pmod2,
\end{equation}
\begin{equation}\label{f2}
 a_{2,7}\Big(2\cdot p^{2\alpha+1}(pn+j)+\dfrac{5\cdot p^{2\alpha+2}-1}{4}\Big)\equiv 0\pmod2,
\end{equation}
\begin{equation}\label{f1a}
\sum _{n=0}^{\infty} a_{2,7}\Big(14\cdot p^{2\alpha}n+\dfrac{21\cdot p^{2\alpha}-1}{4}\Big)q^n\equiv f_2f_7\pmod2
\end{equation} and
\begin{equation}\label{f2a}
 a_{2,7}\Big(14\cdot p^{2\alpha+1}(pn+j)+\dfrac{21\cdot p^{2\alpha+2}-1}{4}\Big)\equiv 0\pmod2.
\end{equation}
\end{theorem}

\begin{theorem}\label{thm7} 
If $w\in{ \left\lbrace 13, 17\right\rbrace}$, then for all $\alpha\geq0,$ we have
\begin{equation}\label{c3}
\sum _{n=0}^{\infty} a_{4,5}\Big(2\cdot 5^{\alpha}n+\dfrac{5^{\alpha}-1}{2}\Big)q^n\equiv f_1f_5 \pmod2
\end{equation} and
\begin{equation}\label{c4}
 a_{4,5}\Big(2\cdot 5^{\alpha+1}n+\dfrac{w\cdot 5^{\alpha}-1}{2}\Big)\equiv 0\pmod2.
\end{equation}
\end{theorem}

\begin{theorem}\label{thm8} Let $p>5$ be a prime with $\Big(\dfrac{-5}{p}\Big)=-1$ and $1\leq j\leq p-1$. Then for  all $\alpha\ge 0$, we have
\begin{equation}\label{c1}
\sum _{n=0}^{\infty} a_{4,5}\Big(2\cdot p^{2\alpha}n+\dfrac{p^{2\alpha}-1}{2}\Big)q^n\equiv f_1f_5\pmod2
\end{equation} and
\begin{equation}\label{c2}
 a_{4,5}\Big(2\cdot p^{2\alpha+1}(pn+j)+\dfrac{p^{2\alpha+2}-1}{2}\Big)\equiv0\pmod2.
\end{equation}
\end{theorem}

\begin{theorem}\label{thm9}Let $w_1\in{ \left\lbrace 3, 5\right\rbrace}$, $w_2\in{ \left\lbrace 13,25,37\right\rbrace}$ and $\alpha\geq0$. Then
\begin{equation}\label{i1}
 a_{4,9}\Big(6n+w_1\Big)\equiv 0\pmod2,
\end{equation}
\begin{equation}\label{i2}
 a_{4,9}\Big(24n+19\Big)\equiv 0\pmod2,
\end{equation}
\begin{equation}\label{i3}
 a_{4,9}\Big(6\cdot 4^{\alpha+2}n+20\cdot 4^{\alpha+1}-1\Big)\equiv 0\pmod2,
\end{equation}
\begin{equation}\label{i4}
 a_{4,9}\Big(48n+w_2\Big)\equiv 0\pmod2
\end{equation} and 
\begin{equation}\label{i5}
a_{4,9}(48n+1)
\equiv \left\{
        \begin{array}{ll}
         1 \pmod2   
           & if~ n~ is~ a~ pentagonal~ number, \\
           0 \pmod2
           & otherwise. 
        \end{array}
    \right.
\end{equation}
\end{theorem}

\section{Preliminaries} 
 In this section, we collect the $q$-series identities that are used in our proofs. Recall that Ramanujan's general theta-function $f\left(a,b\right)$ is defined by
 \begin{equation}\label{eq2}
  f(a, b)=\sum_{n=-\infty}^\infty a^{n(n+1)/2}b^{n(n-1)/2}.
 \end{equation} 
 Important special cases of $f(a, b)$ \cite[p. 36, Entry 22 (i), (ii), (iii)]{bc3} are  the theta-functions $\phi(q)$, $\psi(q)$ and $f(-q)$, which satisfy the identities
 \begin{equation}\label{e0}
\phi(q):=f(q, q)=\sum_{n=0}^{\infty}q^{n^2}=(-q;q^2)^{2}_{\infty}(q^2;q^2)_{\infty}=\dfrac{{f
_2}^5}{{f_1}^2{f_4}^2}, 
 \end{equation} 
\begin{equation}\label{l1}
\psi(q):=f(q,q^3)=\sum_{n=0}^{\infty}q^{n(n+1)/{2}}=\dfrac{(q^2;q^2)_{\infty}}{(q;q^2)_{\infty}}=\dfrac{{f
_2}^2}{f_1},
\end{equation} and
\begin{equation}\label{l2} 
 f(-q):=f(-q, -q^2)=\sum_{n=0}^{\infty}(-1)^n q^{n(3n-1)/{2}}=(q;q)_{\infty}=f_1 .\end{equation}
 In terms of $f(a, b)$, Jacobi's triple product identity \cite[Entry 19, p.35]{bc3}   is given by
 \begin{equation}\label{h20}
 f(a,b)=(-a;ab)_{\infty}(-b;ab)_{\infty} (ab;ab)_{\infty}.\end{equation}

\begin{lemma}
   \cite[Theorem 2.2]{cg} For any prime $p\geq 5$, we have
  \begin{multline}\label{u10} 
f_1=\sum_{\substack{ k={-(p-1)/2} \\ k \ne {(p^\ast-1)/6}}}^{k={(p-1)/2}}(-1)^{k}{\it q}^{(3k^2+k)/2} f\left( {-\it q}^{{(3p^2+(6k+1)p)/2}},{-\it q}^{(3p^2-(6k+1)p)/2}\right) \\+(-1)^{{(p^\ast-1)/6}}{\it q}^{(p^2-1)/24}f_{p^2},\end{multline}where
\begin{equation*}
 p^\ast
= \left\{
        \begin{array}{ll}
         p,   
           & if~ p \equiv 1\pmod 6 \\
           -p ,
           & if~ p \equiv 5\pmod 6. 
        \end{array}
    \right.
\end{equation*} \\Furthermore,  if $$\dfrac{-(p-1)}{2}\leq k \leq\dfrac{(p-1)}{2}~and ~k \neq \dfrac{(p^\ast-1)}{6},$$~ then\\
$$\dfrac{3k^2+k}{2}\not\equiv \dfrac{p^2-1}{24} \pmod p.$$
\end{lemma}

\begin{lemma}\cite[p. 303, Entry 17(v)]{bc3}  We have that
\begin{equation}\label{u7}
f_1=f_{49}\left(\dfrac{B(q^7)}{C(q^7)}-q \dfrac{A(q^7)}{B(q^7)}-q^2+q^5\dfrac{C(q^7)}{A(q^7)}\right),
\end{equation}
where $A(q)=f(-q^3,-q^4),  B(q)=f(-q^2,-q^5)~and ~C(q)=f(-q,-q^6)$.
\end{lemma}

\begin{lemma} \cite{MD} We have that \begin{equation} \label{g1}
f_1=f_{25}(R(q^5)-q-q^2R(q^5)^{-1}),
\end{equation} where $$R(q)=\dfrac{(q^2;q^5)_{\infty}(q^3;q^5)_{\infty}}{(q;q^5)_{\infty}(q^4;q^5)_{\infty}}.$$
\end{lemma}

\begin{lemma}\label{i}We have
\begin{equation}\label{t1}
\dfrac{1}{{f_1^2}}=\dfrac{{f_8^5}}{{f_2^5}{f_{16}^2}}+2q\dfrac{{f_4^2}{f_{16}^{2}}}{{f_2^5}{f_8}},\end{equation}
\begin{equation}\label{t2}
\dfrac{f_9}{f_1}= \dfrac{{f_{12}^3}{f_{18}}}{{f_2^2}{f_6}{f_{36}}}+q \dfrac{{{f_4^2}}{f_6}{f_{36}}}{{f_2^3}{{f_{12}}}},
\end{equation}
\begin{equation}\label{t3}
f_{1}f_{5}^3=f_2^3f_{10}-q\dfrac{f_2^2f_{10}^2f_{20}}{f_4}+2q^2f_4f_{20}^3-2q^3\dfrac{f_4^4f_{10}f_{40}^2}{f_2f_8^2},
\end{equation}
\begin{equation}\label{t4}
\dfrac{f_5}{f_1}= \dfrac{{f_8}{f_{20}^2}}{f_2^2f_{40}}+q \dfrac{f_4^3f_{10}f_{40}}{f_2^3 f_8f_{20}},
\end{equation}
\begin{equation}\label{t5}
\dfrac{f_3^3}{f_1}= \dfrac{{f_4^3 f_6^2}}{f_2^2 f_{12}}+q \dfrac{f_{12}^3}{f_4},
\end{equation}
\begin{equation}\label{t6}
f_{1}f_{7}=\dfrac{f_2 f_{14}f_{16}^2f_{56}^5}{f_4f_8f_{28}^3f_{112}^2}-qf_4f_{28}+q^6\dfrac{f_2f_8^5f_{14}f_{112}^2}{f_4^3 f_{16}^2f_{28}f_{56}}.
\end{equation}
\end{lemma}

For the proof of \eqref{t1}, see Hirschhorn \cite[p.40] {HS4}. Equation  \eqref{t2} was proved by Xia and Yao \cite{XY}. For  the proof  of \eqref{t3}, see Naika et.al \cite{MBH}. Equation \eqref{t4} was proved by Hirschhorn and Sellers \cite {HS}. Equation \eqref{t5} was proved by Hirschhorn et.al \cite {HGB}. Equation \eqref{t6} was proved by Xia \cite[Lemma 3.14]{xia5}.

To end this section, we record the following congruence which can be easily proved using the binomial theorem: For all positive integers $t$ and $m$ we have 
\begin{equation}\label{y7}
{ f_t}^{2m}\equiv {f_{2t}}^m\pmod2.\end{equation}

\section{\bf Proof of Theorems \ref{thm1}-\ref{thm4} }
\noindent\textbf{Proof of Theorem \ref{thm1}:}\begin{proof}
Setting $(r,s)=(2,5)$ in \eqref{eq1} and using elementary $q$-operations, we obtain
\begin{equation}\label{s1}
\sum _{n=0}^{\infty}a_{2,5}(n)q^n=\dfrac{f_2^2f_5f_{20}}{f_1f_4f_{10}^2}.
\end{equation}Combining \eqref{t4} and \eqref{s1}, we find that
\begin{equation}\label{j1}
\sum _{n=0}^{\infty}a_{2,5}(n)q^n=\dfrac{f_8f_{20}^3}{f_4f_{10}^2f_{40}}+q\dfrac{f_4^2f_{40}}{f_2f_8f_{10}}.\end{equation}
 Extracting the terms involving even powers of $q$ of \eqref{j1}, we obtain
\begin{equation}\label{s2}
\sum _{n=0}^{\infty}a_{2,5}(2n)q^n=\dfrac{f_4f_{10}^3}{f_2f_5^2f_{20}}.
\end{equation} In view of \eqref{y7}, \eqref{s2} can be written as 
\begin{equation}\label{s3}
\sum _{n=0}^{\infty}a_{2,5}(2n)q^n\equiv f_2 \pmod2.
\end{equation}Extracting the terms involving odd powers of $q$ from \eqref{s3} yields \eqref{e1}. Finally, extracting the terms involving even powers of $q$ from both sides of \eqref{s3} and using \eqref{l2} yields \eqref{e2}.
\end{proof}

\noindent \textbf{Proof of Theorem \ref{thm3}:} \begin{proof}
Extracting the terms involving odd powers of $q$ from both sides of \eqref{j1}, we obtain
\begin{equation}\label{s4}
\sum _{n=0}^{\infty}a_{2,5}(2n+1)q^n=\dfrac{f_2^2f_{20}}{f_1f_4f_5}.\end{equation} 
In view of \eqref{y7}, we can rewrite \eqref{s4} as
\begin{equation}\label{s5}
\sum _{n=0}^{\infty}a_{2,5}(2n+1)q^n\equiv\dfrac{f_1f_5^3}{f_2} \pmod2.
\end{equation}
 Combining \eqref{t3} and \eqref{s5}, we find that
\begin{equation}\label{s6}
\sum _{n=0}^{\infty}a_{2,5}(2n+1)q^n\equiv f_2^2f_{10}-q\dfrac{f_2f_{10}^2f_{20}}{f_4} \pmod2.
\end{equation}
Extracting the terms involving even powers of $q$ from both sides of \eqref{s6}, we obtain
\begin{equation}\label{s7}
\sum _{n=0}^{\infty}a_{2,5}(4n+1)q^n\equiv f_2f_5 \pmod2.
\end{equation}
Equation \eqref{s7} is the $\gamma=0$ case of \eqref{e8}. Now suppose that \eqref{e8} holds for some $\gamma\geq0$. Using \eqref{u10} in \eqref{e8}, we deduce that 
\begin{multline}\label{i13}
 \sum_{n\geq0}^{\infty} a_{2,5}\Big(4\cdot p^{2\gamma}n+\dfrac{7\cdot p^{2\gamma}-1}{6}\Big) q^n \\ \equiv  \Bigg[ \sum_{\substack{ k={-(p-1)/2} \\ k \ne {(p^\ast-1)/6}}}^{k={(p-1)/2}}{\it q}^{3k^2+k} f\big( {-\it q}^{3p^2+(6k+1)p},{-\it q}^{3p^2-(6k+1)p}\big ) \\ 
 +{\it q}^{(p^2-1)/12}f_{2p^2} \Bigg]\\ \times  \Bigg[ \sum_{\substack{ m={-(p-1)/2} \\ m \ne {(p^\ast-1)/6}}}^{m={(p-1)/2}}{\it q}^{5(3m^2+m)/2} f\left( {-\it q}^{5(3p^2+(6m+1)p)/2},{-\it q}^{5(3p^2-(6m+1)p)/2}\right)\\   
 +{\it q}^{5(p^2-1)/24}f_{5 p^2} \Bigg] \pmod 2. \end{multline} 
 Consider the congruence
$$3k^2+k+5\dfrac{(3m^2+m)}{2} \equiv \dfrac{7(p^2-1)}{24} \pmod p,$$
which is equivalent to 
$$ (12k+2)^2+10(6m+1)^2 \equiv 0 \pmod p.$$
Since $\Big(\dfrac{-10}{p}\Big)=-1$, the only solution of this congruence is $k=m=\dfrac{(p^\ast-1)}{6}$. Therefore, extracting the terms involving $q^{pn+{{7(p^2-1)/24}}}$ from both sides of
\eqref{i13}, dividing by $q^{{7(p^2-1)/24}}$ and then replacing $q^{p}$ by $q$, we find that
\begin{equation}\label{s12}
\sum _{n=0}^{\infty} a_{2,5}\Big(4\cdot p^{2\gamma+1}n+\dfrac{7\cdot p^{2\gamma+2}-1}{6}\Big)q^n\equiv f_{2p}f_{5p}\pmod2,
\end{equation}
  which yields
 \begin{equation}\label{s13}
\sum _{n=0}^{\infty} a_{2,5}\Big(4\cdot p^{2(\gamma+1)}n+\dfrac{7\cdot p^{2(\gamma+1)}-1}{6}\Big)q^n\equiv f_2f_5 \pmod2,
\end{equation} 
which is the $\gamma +1$ case of \eqref{e8}. On the other hand, extracting the terms involving $q^{pn+j}~(1\leq j \leq p-1)$ from \eqref{s12}, we arrive at \eqref{e9}.
Employing \eqref{g1} in \eqref{e8} and then extracting the  terms involving $q^{5n+2}$ yields \eqref{e8a}. Next, using \eqref{u10} in \eqref{e8a} and proceeding as in the proof of \eqref{e8}, we arrive at \begin{equation}\label{yu1}
\sum _{n=0}^{\infty} a_{2,5}\Big(20\cdot p^{2\gamma+1}n+\dfrac{55\cdot p^{2\gamma+2}-1}{6}\Big)q^n\equiv f_pf_{10p} \pmod2.
\end{equation}
Finally, \eqref{e9b} follows from extracting the terms involving $q^{pn+j}~(1\leq j \leq p-1)$ from \eqref{yu1}. 
\end{proof}

\noindent\textbf{Proof of Theorem \ref{thm4}:} 
\begin{proof}
Combining \eqref{t1} and \eqref{s2}, we find that
\begin{equation}\label{s14}
\sum _{n=0}^{\infty}a_{2,5}(2n)q^n=\dfrac{f_4f_{40}^5}{f_2f_{10}^2f_{20}f_{80}^2}+2q^5\dfrac{f_4f_{20}f_{80}^2}{f_2f_{10}^2f_{40}}.\end{equation}
Extracting the terms involving odd powers of $q$ from both sides of \eqref{s14}, we obtain
\begin{equation}\label{s15}
\sum _{n=0}^{\infty}a_{2,5}(4n+2)q^n=2q^2\dfrac{f_2f_{10}f_{40}^2}{f_1f_{5}^2f_{20}}.\end{equation}
In view of \eqref{y7}, \eqref{s15} can be written as
\begin{equation}\label{s16}
\sum _{n=0}^{\infty}a_{2,5}(4n+2)q^n\equiv2q^2f_1f_{20}^3 \pmod4,\end{equation}
 which is the $\beta=0$ case  of \eqref{e10}. Now assume that \eqref{e10} holds for some  $\beta\geq0$. Employing \eqref{g1} in \eqref{e10}, we arrive at  
\begin{equation}\label{s17}
\sum _{n=0}^{\infty} a_{2,5}\Big(4\cdot 5^{2\beta}n+\dfrac{13\cdot 5^{2\beta}-1}{6}\Big)q^n\equiv2q^2 f_{20}^3 f_{25}\Big(R(q^{5})-q-q^2R(q^{5})^{-1})\Big)\pmod4.
\end{equation}
 Extracting the terms involving $q^{5n+3}$ from \eqref{s17} yields \eqref{e11}. Next, using \eqref{g1} in \eqref{e11} and then extracting the terms involving $q^{5n+2}$, we obtain  
\begin{equation}\label{s19}
\sum _{n=0}^{\infty} a_{2,5}\Big(4\cdot 5^{2(\beta+1)}n+\dfrac{13\cdot 5^{2(\beta+1)}-1}{6}\Big)q^n\equiv2q^2 f_1f_{20}^3 \pmod4,
\end{equation}
which is the $\beta +1$ case of \eqref{e10}. 
Employing \eqref{g1} in \eqref{e10} and then
extracting the terms involving $q^{5n+j}$ for $j\in{\{0,1\}}$  yields \eqref{e12}. Finally, using \eqref{g1} in \eqref{e11} and then extracting terms involving $q^{5n+j}$ for $j \in {\{1,3\}}$ yields \eqref{e13}.
 \end{proof}
 
\section{\bf Proof of Theorem \ref{thm5} }
\begin{proof}
Setting $(r,s)=(2,7)$ in \eqref{eq1} and using elementary $q$-operations, we have
\begin{equation}\label{s20}
\sum _{n=0}^{\infty}a_{2,7}(n)q^n=\dfrac{f_2^2f_7f_{28}}{f_1f_4f_{14}^2}.
\end{equation} 
In view of \eqref{y7}, we can rewrite \eqref{s20} as
\begin{equation}\label{s21}
\sum _{n=0}^{\infty}a_{2,7}(n)q^n\equiv\dfrac{f_1f_7}{f_2}\pmod2.\end{equation}
Combining \eqref{t6} and \eqref{s21}, we find that\begin{equation}\label{s22}
\sum _{n=0}^{\infty}a_{2,7}(n)q^n\equiv\dfrac{f_{14}f_{16}^2f_{56}^5}{f_4f_8f_{28}^3f_{112}^2}-q\dfrac{f_4f_{28}}{f_2}+q^6\dfrac{f_8^5f_{14}f_{112}^2}{f_4^3 f_{16}^2f_{28}f_{56}}\pmod2.\end{equation}
Extracting the terms involving odd powers of $q$ from both sides of \eqref{s22}, we obtain
\begin{equation}\label{s23}
\sum _{n=0}^{\infty}a_{2,7}(2n+1)q^n\equiv{f_1f_{14}}\pmod2,\end{equation}
which is the $\alpha=0$ case of \eqref{f1}. Using \eqref{u10} in \eqref{f1} and proceeding as in the proof of  \eqref{e8}, we arrive at 
\begin{equation}\label{s27}
\sum _{n=0}^{\infty} a_{2,7}\Big(2\cdot p^{2\alpha+1}n+\dfrac{5\cdot p^{2\alpha+2}-1}{4}\Big)q^n\equiv f_{p}f_{14p}\pmod2,
\end{equation} 
which yields
\begin{equation}\label{s28}
\sum _{n=0}^{\infty} a_{2,7}\Big(2\cdot p^{2(\alpha+1)}n+\dfrac{5\cdot p^{2(\alpha+1)}-1}{4}\Big)q^n\equiv f_1f_{14}\pmod2,
\end{equation}
which is the $\alpha +1$ case of \eqref{f1}. On the other hand, extracting the terms involving $q^{pn+j}~(1\leq j \leq p-1)$ from \eqref{s27}, we arrive at \eqref{f2}. Next,
using \eqref{u7} in \eqref{f1} and then extracting the terms involving $q^{7n+2}$ yields \eqref{f1a}. Now 
employing \eqref{u10} in \eqref{f1a} and proceeding as in the proof of  \eqref{f1}, we arrive at \begin{equation}\label{jk1}
\sum _{n=0}^{\infty} a_{2,7}\Big(14\cdot p^{2\gamma+1}n+\dfrac{21\cdot p^{2\gamma+2}-1}{4}\Big)q^n\equiv f_{2p}f_{7p} \pmod2.
\end{equation}
Finally, \eqref{f2a} follows from extracting the terms involving $q^{pn+j} (1\leq j \leq p-1)$ from \eqref{jk1}.
 \end{proof} 
 
\section{\bf Proof of Theorems \ref{thm7}-\ref{thm8} }
\noindent\textbf{Proof of Theorem \ref{thm7}:}
\begin{proof}
Setting $(r,s)=(4,5)$ in \eqref{eq1} and using elementary $q$-operations, we obtain
\begin{equation}\label{s29}
\sum _{n=0}^{\infty}a_{4,5}(n)q^n=\dfrac{f_2f_4f_5f_{40}}{f_1f_8f_{10}f_{20}}.
\end{equation}
Combining \eqref{t4} and \eqref{s29}, we find that
\begin{equation}\label{s30}
\sum _{n=0}^{\infty}a_{4,5}(n)q^n=\dfrac{f_4f_{20}}{f_2f_{10}}+q\dfrac{f_4^4f_{40}^2}{f_2^2f_8^2f_{20}^2}.
\end{equation}
Extracting the terms involving even powers of $q$ from both sides of \eqref{s30}, we obtain
\begin{equation}\label{s31}
\sum _{n=0}^{\infty}a_{4,5}(2n)q^n=\dfrac{f_2f_{10}}{f_1f_5}.
\end{equation}
In view of \eqref{y7}, we can rewrite \eqref{s31} as
\begin{equation}\label{s32}
\sum _{n=0}^{\infty}a_{4,5}(2n)q^n\equiv f_1f_5\pmod2,
\end{equation}
which is the $\alpha=0$ case of \eqref{c3}. Now assume that \eqref{c3} holds for some  $\alpha\geq0$. Using \eqref{g1} in \eqref{c3}, we find that 
\begin{equation}\label{s33a}
\sum _{n=0}^{\infty} a_{4,5}\Big(2\cdot 5^{\alpha}n+\dfrac{5^{\alpha}-1}{2}\Big)q^n\equiv f_5 f_{25}\Big(R(q^{5})-q-q^2R(q^{5})^{-1})\Big)\pmod2.
\end{equation} 
Extracting the terms involving $q^{5n+1}$ from both sides of \eqref{s33a}, we arrive at
\begin{equation}\label{s33}
\sum _{n=0}^{\infty} a_{4,5}\Big(2\cdot 5^{\alpha+1}n+\dfrac{5^{\alpha+1}-1}{2}\Big)q^n\equiv f_1f_5 \pmod2,
\end{equation}
which is the $\alpha +1$ case of \eqref{c3}. Finally, using \eqref{g1} in \eqref{c3} and then extracting the terms involving $q^{5n+j}$ for $j\in {\{3,4\}}$ yields \eqref{c4}.
\end{proof}

\noindent\textbf{Proof of Theorem \ref{thm8}:}
\begin{proof}
 Congruence \eqref{s32} is the $\alpha=0$ case of \eqref{c1}. Now suppose that \eqref{c1} holds for some $\alpha\geq0$. Using \eqref{u10} in \eqref{c1} and proceeding as in the proof of  \eqref{e8}, we arrive at 
 \begin{equation}\label{s34}
 \sum _{n=0}^{\infty} a_{4,5}\Big(2\cdot p^{2\alpha+1}n+\dfrac{p^{2\alpha+2}-1}{2}\Big)q^n\equiv f_{p}f_{5p}\pmod2,
 \end{equation}
 which yields 
 \begin{equation}\label{s35}
 \sum _{n=0}^{\infty} a_{4,5}\Big(2\cdot p^{2(\alpha+1)}n+\dfrac{p^{2(\alpha+1)}-1}{2}\Big)q^n\equiv f_1f_{5}\pmod2,
\end{equation}
which is the case $\alpha +1$  of \eqref{c1}. On the other hand, extracting the terms involving $q^{pn+j} (1\leq j\leq p-1)$ from \eqref{s34}, we arrive at \eqref{c2}.
\end{proof}

\section{\bf Proof of Theorem \ref{thm9} }
\begin{proof}Setting $(r,s)=(4,9)$ in \eqref{eq1} and using elementary $q$-operations, we obtain
\begin{equation}\label{s36}
\sum _{n=0}^{\infty}a_{4,9}(n)q^n=\dfrac{f_2f_4f_9f_{72}}{f_1f_8f_{18}f_{36}}.
\end{equation}
Combining \eqref{t2} and \eqref{s36}, we find that
\begin{equation}\label{s37}
\sum _{n=0}^{\infty}a_{4,9}(n)q^n=\dfrac{f_4f_{12}^3f_{72}}{f_2f_6f_8f_{36}^2}+ q\dfrac{f_4^3f_6f_{72}}{f_2^2f_8f_{12}f_{18}}.
\end{equation}
Extracting  the terms involving odd powers of $q$ from \eqref{s37} and then employing \eqref{y7}, we obtain
\begin{equation}\label{s38}
\sum _{n=0}^{\infty}a_{4,9}(2n+1)q^n\equiv\dfrac{f_9^3}{f_3}\pmod2.
\end{equation}
 Comparing the  terms involving $q^{3n+j}$, for $j\in{\{1,2\}}$ from both sides of \eqref{s38} yields \eqref{i1}. Next, extracting the  terms involving $q^{3n}$ from  both sides of \eqref{s38}, we obtain
\begin{equation}\label{s39}
\sum _{n=0}^{\infty}a_{4,9}(6n+1)q^n\equiv\dfrac{f_3^3}{f_1}\pmod2.
\end{equation} 
Combining \eqref{t5} and \eqref{s39}, we find that
\begin{equation}\label{s40}
\sum _{n=0}^{\infty}a_{4,9}(6n+1)q^n\equiv\dfrac{f_4^3f_6^2}{f_2^2f_{12}}+q\dfrac{f_{12}^3}{f_4}\pmod2.
\end{equation}
Extracting the terms involving odd powers of $q$ from \eqref{s40}, we obtain
\begin{equation}\label{s41}
\sum _{n=0}^{\infty}a_{4,9}(12n+7)q^n\equiv\dfrac{f_6^3}{f_2}\pmod2.
\end{equation}
 Extracting the terms involving odd powers of $q$ from  \eqref{s41} yields \eqref{i2}. Next, extracting the terms involving even powers of $q$ from both sides of \eqref{s41}, we obtain
\begin{equation}\label{s42}
\sum _{n=0}^{\infty}a_{4,9}(24n+7)q^n\equiv\dfrac{f_3^3}{f_1}\pmod2.
\end{equation}
Combining \eqref{s39} and \eqref{s42}, we find that 
\begin{equation}\label{s43}
a_{4,9}(24n+7)\equiv a_{4,9}(6n+1)\pmod2.
\end{equation}
 From \eqref{s43} and by mathematical induction, we have \begin{equation}\label{s44}
 a_{4,9}\Big(6\cdot 4^{\alpha+1}n+2\cdot 4^{\alpha+1}-1\Big)\equiv a_{4,9}(6n+1)\pmod2.
\end{equation}
 Using \eqref{s44} and congruence \eqref{i2}, we arrive at \eqref{i3}. On the other hand, extracting the terms involving even powers of $q$ from both sides of \eqref{s40}, we obtain
\begin{equation}\label{s45}
\sum _{n=0}^{\infty}a_{4,9}(12n+1)q^n \equiv f_4\pmod2.
\end{equation}
Extracting the terms involving $q^{4n+j}$ for $j\in{\{1,2,3 \}}$ from \eqref{s45} yields \eqref{i4}. On the other hand, extracting the terms involving $q^{4n}$ from both sides of \eqref{s45} and using \eqref{l2} yields \eqref{i5}.
\end{proof}
\section*{Acknowledgement} The authors are thankful to the referee for his/her comments which improves the quality of the paper. The first author thanks University Grants Commission (UGC) of India for supporting her research work through Junior Research Fellowship (JRF) vide  UGC-Ref.no. 1107/CSIR-UGC NET DEC-2018.

\end{document}